\input amstex
\documentstyle{amsppt}
\hsize160mm
 \vsize220mm
 \loadbold

\centerline{\bf  On a Hierarchy of Means}

\vskip 1cm

\centerline{ \bf Slavko Simic}
\bigskip
\centerline{ Mathematical Institute SANU}

\centerline {Kneza Mihaila 36, 11000 Belgrade, Serbia}
\bigskip
\centerline{\sl e-mail: ssimic\@turing.mi.sanu.ac.rs}

 \vskip 1cm


\vskip 1cm

{\sevenrm{\bf Abstract} \ For a class of partially ordered means
we introduce a notion of the (nontrivial) cancelling mean. A
simple method is given which helps to determine cancelling means
for well known classes of H\"{o}lder and Stolarsky means.}

\vskip 1cm

{\bf 1. Introduction}

\bigskip

A {\it mean} is a map $M: \Bbb R_+\times\Bbb R_+\to \Bbb R_+$,
with a property

$$
\min(a,b)\le M(a,b)\le\max(a,b),
$$

for each $a,b\in \Bbb R_+$.

\bigskip

Denote by $\Omega$ the class of means which are {\it symmetric}
(in variables $a,b$), {\it reflexive} and {\it homogeneous}
(necessarily of order one). We shall consider in the sequel only
means from this class.

\bigskip

The set of means can be equipped with a partial ordering defined
by $M\le N$ if and only if $M(a,b)\le N(a,b)$ for all $a,b\in \Bbb
R_+$. Thus, $\Delta$ is an {\it ordered} family of means if for
any $M, N\in {\Delta}$ we have $M\le N$ or $N\le M$.

\bigskip

Most known ordered family of means is the following family
$\Delta_0$ of elementary means,

$$
\Delta_0: \  H\le G \le L\le I\le A\le S,
$$
where
$$
H=H(a, b):=2(1/a+1/b)^{-1}; \ \ G=G(a, b):=\sqrt {ab}; \ \ L=L(a,
b):={b-a\over \log b-\log a};
$$
$$
 I=I(a, b):=(b^b/a^a)^{1/(b-a)}/e; \ \ A=A(a, b):={a+b\over 2}; \ \
 S=S(a,b):= a^{a\over a+b} b^{b\over a+b},
$$
are the harmonic, geometric, logarithmic, identric, arithmetic and
Gini  mean, respectively.

\bigskip

Another example is the class of H\"{o}lder (or Power) means
$\{A_s\}$, defined for $s\in \Bbb R$ as

$$
A_s(a, b):=\Bigl({a^s+b^s\over 2}\Bigr)^{1/s}, \ A_0=G.
$$

It is well known that the inequality $A_s(a,b)<A_t(a,b)$ holds for
all $a,b\in \Bbb R_+, a\neq b$ if and only if $s<t$. This property
is used in a number of papers for approximation of a particular
mean by means from the class $\{A_s\}$.
\bigskip
Hence (cf [3], [4], [10]),

$$
G=A_0<L<A_{1/3}; \ \ A_{2/3}<I<A_1=A; \ \
A_{\log_{\pi}2}<P<A_{2/3}; \ \ A_{\log_{\pi/2}2}<T<A_{5/3},
$$

where all bounds are best possible and Seiffert means $P$ and $T$
are defined by

$$
P=P(a,b):={{a-b}\over{2\arcsin{{a-b}\over{a+b}}}}; \ \
T=T(a,b):={{a-b}\over{2\arctan{{a-b}\over{a+b}}}}.
$$
\bigskip
In the recent paper [8] we introduce a more complex structured
class of means $\{\lambda_s\}$, given by

$$
\lambda_s(a,b):={s-1\over s+1}{A_{s+1}^{s+1}-A^{s+1}\over
A_s^s-A^s}, \ s\in\Bbb R,
$$

that is,

$$
\lambda_s(a, b):=\cases {s-1\over s+1}{a^{s+1}+b^{s+1}-2({a+b\over
2})^{s+1}\over a^s+b^s-2({a+b\over 2})^s},& s\in \Bbb
R/\{-1, 0, 1\};\\
                {2\log {a+b\over 2}-\log a-\log b\over {1\over 2a}+{1\over 2b}-{2\over a+b}}, & s=-1;\\
                {a\log a+b\log b-(a+b)\log{a+b\over 2}\over 2\log{a+b\over 2}-\log a-\log b},& s=0;\\
                {(b-a)^2\over 4(a\log a+b\log b-(a+b)\log{a+b\over 2})},& s=1.
                \endcases
                $$
\bigskip

Those means are obviously symmetric and homogeneous of order one.

We also proved that $\lambda_s$ is monotone increasing in
$s\in\Bbb R$; therefore $\{\lambda_s\}$ represents an ordered
family of means.
\bigskip
Among others, the following approximations are obtained for $a\neq
b$:

$$
\lambda_{-4}<H<\lambda_{-3}; \ \lambda_{-1}<G<\lambda_{-1/2}; \
\lambda_{0}<L<\lambda_{1}<I<\lambda_{2}=A; \ \lambda_{5}<S,
$$

and {\it there is no finite $s>5$ such that the inequality
$S(a,b)\le \lambda_s(a,b)$ holds for each $a,b\in \Bbb R^+$}.

\bigskip

This last result shows that, in a sense, the mean $S$ is "greater"
than any other mean from the class $\{\lambda_s\}$. We shall say
that $S$ is the {\it cancelling mean} for the class
$\{\lambda_s\}$.

\bigskip

{\bf Definition 1} \ {\it The mean $S^*(\Delta)$ is right
cancelling mean for an ordered class of means
$\Delta\subset\Omega$ if there exists $M\in \Delta$ such that
$S^*(a,b)\ge M(a,b)$ but there is no mean $N\in \Delta$ such that
the inequality $N(a,b)\ge S^*(a,b)$ holds for each $a,b\in\Bbb
R_+$.}
\bigskip
Definition of the left cancelling mean $S_*$ is analogues.
\bigskip
Therefore,
$$
S_*(\Delta_0)=H; \ \  S^*(\Delta_0)=S; \ \ S^*(\lambda_s)=S.
$$
\bigskip
Of course that the left and right cancelling means exist for
arbitrary ordered family of means as $S^*(a,b)=\max(a,b), \
S_*(a,b)=\min(a,b)$. We call them {\it trivial}.
\bigskip
The aim of this article is to determine non-trivial cancelling
means for some well known classes of ordered means. We shall also
give a simple criteria for the right cancelling mean with further
discussion in the sequel.

 \bigskip
 As an illustration of problems and methods which shall be treated
 in this paper, we prove firstly the following,
 \bigskip
 {\bf 1.1 Cancellation theorem for the Generalized Logarithmic Means}
\bigskip
The family of Generalized Logarithmic Means $\{L_p\}$ is given by

$$
L_p=L_p(a,b):=\Bigl({a^p-b^p\over p(\log a-\log b)}\Bigr)^{1/p}, \
p\in\Bbb R; \ L_0=G, \ L_1=L.
$$

It is a subclass of well-known Stolarsky means (cf [2],[5],[7])
hence symmetric, homogeneous and monotone increasing in $p$.
Therefore it represents an ordered family of means.
\bigskip
{\bf Theorem 1.1} \ {\it For the class $\{L_p\}$ we have

$$
S_*(L_p)=H, \ S^*(L_p)=A.
$$

Moreover, for $-3<p<3, \ a\neq b$,

$$
S_*(L_p)=H(a,b)<L_{-3}(a,b)<L_p(a,b)<L_3(a,b)<A(a,b)= S^*(L_p),
$$

with those bounds as best possible. }
\bigskip
{\bf Proof} \ We prove firstly that the inequality
$L_3(a,b)<A(a,b)$ holds for all $a,b\in \Bbb R_+, \ a\neq b$.
\bigskip
Indeed,
$$
{L_3^3\over A^3}={({2a\over a+b})^3-({2b\over a+b})^3\over
3(\log{2a\over a+b}-\log{2b\over a+b})}={(1+t)^3-(1-t)^3\over
3(\log(1+t)-\log(1-t))}={3+t^2\over 3(1+t^2/3+t^4/5+\dots)}<1,
$$
where we put $t:={a-b\over a+b}, \ -1<t<1$.
\bigskip

Also,

$$
{L_p^p\over A^p}={(1+t)^p-(1-t)^p\over p(\log(1+t)-\log(1-t))},
$$

and

$$
\lim_{t\to 0}{1\over t^2}\Bigl({L_p^p\over A^p}-1\Bigr)={1\over
6}p(p-3).
$$

Thus $p=3$ is the largest $p$ such that the inequality
$L_p(a,b)\le A(a,b)$ holds for each $a,b\in\Bbb R_+$, since for
$p>3$ and $t$ sufficiently small (i.e., $a$ is sufficiently close
to $b$) we have that $L_p(a,b)> A(a,b)$.
\bigskip
We shall show now that $A$ is the right cancelling mean for the
class $\{L_p\}$.
\bigskip
Indeed, since $\lim_{t\to 1^-}{L_p^p\over A^p}=0$ for fixed $p,
p>3$, we conclude that the inequality $L_p\ge A$ cannot hold.

Hence by Definition 1., $A$ is the right cancelling mean for the
class $\{L_p\}$.
\bigskip
Noting that $ H(a,b)={ab\over A(a,b)}$ and $L_{-p}(a,b)={ab\over
L_p(a,b)}$, we readily get
$$
L_{-p}(a,b)\ge L_{-3}(a,b)\ge H(a,b)=S_*(L_p).
$$

\bigskip

{\bf 2. \ Characteristic number and characteristic function} \ Let
$M=M(a,b)$ be an arbitrary homogeneous and symmetric mean. In
order to facilitate determination of a non-trivial right
cancelling mean, we introduce here a notion of {\it characteristic
number} $\sigma (M)$ as
$$
\sigma(M):=\lim_{a/b\to\infty}{M(a,b)\over
A(a,b)}=M(2,0^+)=M(0^+,2).
$$
\bigskip
Because of homogeneousness, we have
$$
{M(a,b)\over A(a,b)}=M({2a\over a+b}, {2b\over a+b})=M(2{{a\over
b}\over {a\over b}+1}, {2\over {a\over b}+1}),
$$

and the result follows.
 \bigskip
Therefore,
$$
\sigma(H)= \sigma(G)= \sigma(L)=0; \  \sigma(I)=2/e; \
\sigma(A)=1; \ \sigma(S)=2,
$$

and, in general,
$$
0\le \sigma(M)\le 2.
$$
\bigskip

Some simple reasoning gives the next,

 \bigskip
 {\bf Theorem 2.0} \ {\it Let $M,N\in\Omega$. If $M\le N$ then $\sigma(M)\le \sigma(N)$ but if $\sigma(M)> \sigma(N)$ then
 the inequality $M\le N$ cannot hold, at least when $a/b$ is sufficiently
 large}.
 \bigskip
 This assertion is especially important in applications.
 \bigskip
 Also,
 $$
{M(a,b)\over A(a,b)}=M({2a\over a+b}, {2b\over a+b})=M(1-{b-a\over
a+b}, 1+{b-a\over a+b})=M(1-t,1+t),
 $$

 where $t:={b-a\over a+b}, \ -1<t<1$.
 \bigskip
 We say that the function $\phi=\phi_M(t):=M(1-t,1+t)$ is {\it
 characteristic function} for $M$ (related to the arithmetic
 mean). If $\phi$ is analytic then, because of $\phi(0)=1, \ \phi(-t)=\phi(t)$, it has power series
 representation of the form
 $$
 \phi(t)=\sum_0^\infty a_n t^{2n}, \ a_0=1, \ 0\le
 t<1.
 $$
 \bigskip
In this way comparison between means reduces to comparison between
their characteristic functions ([8], [10], [11]).
\bigskip
Obviously,

$$
\phi_H(t)=1-t^2; \ \ \phi_G(t)=\sqrt {1-t^2}; \ \
\phi_L(t)={2t\over \log(1+t)-\log(1-t)}; \ \ \phi_A(t)=1; \eqno(2)
$$
$$
\phi_I(t)=\exp({(1+t)\log(1+t)-(1-t)\log(1-t)\over 2t}-1); \
\phi_S(t)=\exp({1\over 2}((1+t)\log(1+t)+(1-t)\log(1-t))).
$$
\bigskip
Note that
$$
\sigma(M)=\lim_{t\to 1^-}\phi_M(t).
$$
\bigskip
We shall give now some applications of the above.
\bigskip
First of all, for an arbitrary mean $M=M(a,b)$ it is not difficult
to show that $M_s=M_s(a,b):=(M(a^s,b^s))^{1/s}$ is also a mean for
$ s\neq 0$. Especially $M_{-1}(a,b)={ab\over M(a,b)}$ is a mean.
\bigskip
Moreover, it is proved in [9] that the condition $[\log
M(x,y)]_{xy}<0$ is sufficient for $M_s$ to be monotone increasing
in $s\in\Bbb R$ and, if $M\in\Omega$ then $M_0=\lim_{s\to 0}
M_s=G$.
\bigskip
For the family of means $\{M_s\}$ we can state the following {\it
cancellation} assertion.
\bigskip
{\bf Theorem 2.1} \ {\it Let $M\in \Omega$ with $[\log
M(x,y)]_{xy}<0$ and $0<\sigma(M)<2$.
\bigskip
For the ordered class of means
$$
M_s=M_s(a,b):=(M(a^s,b^s))^{1/s}\in \Omega, \ s\neq 0; \ M_0=G,
$$
we have
$$
S_*(M_s)=a^{b\over a+b}b^{a\over a+b}; \ S^*(M_s)=a^{a\over
a+b}b^{b\over a+b}.
$$}
\bigskip
{\bf Proof} \ For fixed $s, s>0$, we have $G=M_0\le M_s$.
\bigskip
But,
$$
\sigma(M_s)=(M(0^+,2^s))^{1/s}=2^{1-1/s}(\sigma(M))^{1/s}<2=\sigma(S).
$$

Therefore, by Theorem 2.0 we conclude that $S$ is the right
cancelling mean for $\{M_s\}$.
\bigskip
Also $G=M_0\ge M_{-s}$. Since
$$
M_{-s}(a,b)=(M(a^{-s},b^{-s}))^{-1/s}=(M((ab)^{-s}b^s,(ab)^{-s}a^s))^{-1/s}=ab(M(b^s,a^s))^{-1/s}={ab\over
M_s(a,b)},
$$
and
$$
a^{b\over a+b}b^{a\over a+b}=a^{1-{a\over a+b}}b^{1-{b\over
a+b}}={ab\over S(a,b)},
$$
it easily follows that $a^{b\over a+b}b^{a\over a+b}=S_*(M_s)$.

\vskip 1cm

Another consequence is the {\it cancellation} assertion for the
family of H\"{o}lder means $A_r=A_r(a,b):=
(A(a^r,b^r))^{1/r}=({a^r+b^r\over 2})^{1/r}, \ A_0=G$. Since
$[\log A(x,y)]_{xy}=-{1\over (x+y)^2}<0$, we obtain (as is already
stated) that $A_r$ are monotone increasing with $r$.
\bigskip
{\bf Theorem 2.2} \ {\it \ For $-2\le r\le 2$ we have
  $$
S_*(A_{r})=a^{b\over a+b}b^{a\over a+b}\le A_{-2}(a,b)\le
A_{r}(a,b)\le A_{2}(a,b)\le a^{a\over a+b}b^{b\over
a+b}=S^*(A_{r}),
  $$

  where given constants are best possible.}
  \bigskip
  {\bf Proof} \ We have
$$
{A_r(a,b)\over S(a,b)}={A_r(a,b)/A(a,b)\over
S(a,b)/A(a,b)}={\phi_{A_r}(t)\over \phi_S(t)},
$$
and
$$
f_r(t):=\log {\phi_{A_r}(t)\over \phi_S(t)}={1\over
r}\log\Bigl({(1+t)^r+(1-t)^r\over 2} \Bigr)-{1\over 2}((1+t)\log
(1+t)+(1-t)\log(1-t)), \ 0< t<1.
$$
\bigskip
Denote
$$
g(t):=2f_2(t)=2\log {\phi_{A_2}(t)\over
\phi_S(t)}=\log(1+t^2)-(1+t)\log (1+t)-(1-t)\log(1-t).
$$

Since

$$
g'(t)={2t\over 1+t^2}-\log(1+t)+\log(1-t),
$$

and

$$
g''(t)={2\over 1+t^2}-{4t^2\over (1+t^2)^2}-{1\over 1+t}-{1\over
1-t}=-{8t^2\over (1+t^2)(1-t^4)}<0,
$$

we clearly have $g'(t)<g'(0)=0$ and $g(t)<g(0)=0$.
\bigskip
Therefore, the inequality $A_2(a,b)\le S(a,b)$ holds for all
$a,b\in \Bbb R_+$.
\bigskip
Also, since
$$
\lim_{t\to 0^+}{f_r(t)\over t^2}={1\over 2}(r-2),
$$

 we conclude that $r=2$ is best possible upper bound for $A_r\le
 S$ to hold.
 \bigskip
 Values for $S_*(A_r)$ and $S^*(A_r)$ follow from Theorem 2.1.

 \vskip 1cm

{\bf 3. \ Cancellation theorem for the class of Stolarsky means}

  \bigskip

 There is a plenty of papers (cf [2], [5], [7]) studying
different properties of the so-called Stolarsky (or extended)
two-parametric mean value, defined for positive values of $x, y,
x\neq y$ by the following
$$
I_{r,s}=I_{r,s}(x,y): =\cases \Bigl({r(x^s-y^s)\over s(x^r-y^r)}\Bigr)^{1/(s-r)},& rs(r-s)\neq 0 \\
                    \exp\Bigl(-{1\over s}+{x^s\log x-y^s\log y\over x^s-y^s}\Bigr), & r=s\neq 0 \\
                    \Bigl({x^s-y^s\over s(\log x-\log y)}\Bigr)^{1/s}, & s\neq 0, r=0 \\
                    \sqrt {xy}, & r=s=0, \\
                    x, & y=x>0.
                    \endcases
                    $$
\bigskip

In this form it was introduced by K. Stolarsky in [5].
\bigskip
Most of the classical two variable means are special cases of the
class $\{I_{r,s}\}$. For example, $I_{1,2}=A$ ,
$I_{0,0}=I_{-1,1}=G$, $I_{-2,-1}=H$ , $I_{0,1}=L$, $I_{1, 1}=I$,
etc.
\bigskip
Main properties of Stolarsky means are given in the following
assertion.
\bigskip
{\bf Proposition 3.1} \ {\it Means $I_{r,s}(x,y)$ are
\bigskip
a. \ symmetric in both parameters, i.e.
$I_{r,s}(x,y)=I_{s,r}(x,y)$;
\bigskip
b. \ symmetric in both variables, i.e.
$I_{r,s}(x,y)=I_{r,s}(y,x)$;
\bigskip
c. \ homogeneous of order one, that is $I_{r,s}(tx,ty)=
tI_{r,s}(x,y), \ t>0$;
\bigskip
d. \ monotone increasing in either $r$ or $s$;
\bigskip
e. \ monotone increasing in either $x$ or $y$;
\bigskip
f. \ logarithmically convex for $r, s\in \Bbb R_-$ and
logarithmically concave for $r,s\in \Bbb R_+$.}

\vskip 1cm

  {\bf Theorem 3.2} \ {\it \ For $-3\le r\le s\le 3$ we have
  $$
S_*(I_{r,s})=a^{b\over a+b}b^{a\over a+b}\le I_{-3,-3}(a,b)\le
I_{r,s}(a,b)\le I_{3,3}(a,b)\le a^{a\over a+b}b^{b\over
a+b}=S^*(I_{r,s}),
  $$

  where given constants are best possible.}
  \bigskip
  {\bf Proof} \ We prove firstly that $I_{3,3}(a,b)\le S(a,b)$ and
  that $s=3$ is the largest constant such that the inequality $I_{s,s}(a,b)\le
  S(a,b)$ holds for all $a,b\in \Bbb R_+$. For this aim we need a
  notion of Lehmer means $l_r$ defined by
  $$
l_r=l_r(a,b):={a^{r+1}+b^{r+1}\over a^r+b^r}.
  $$

  They are continuous and strictly increasing in $r\in \Bbb R$ (cf
  [11]).
  \bigskip
  {\bf Lemma 3.3} ([11]) \ {\it \ $L(a,b)>l_{-{1\over 3}}(a,b)$ for all $a,b>0$ with $a\neq b$, and $l_{-{1\over 3}}(a,b)$ is the best
   possible lower Lehmer mean bound for the logarithmic mean $L(a,b)$.}
   \bigskip
   We also need the following interesting identity which is new to our modest
   knowledge.
   \bigskip
   {\bf Lemma 3.4} \ {\it For all $s\in \Bbb R/\{0\}$ we have
   $$
\log{I_{s,s}(a,b)\over S(a,b)}={1\over s}\Bigl({l_{-{1\over
s}}(a^s,b^s)\over L(a^s,b^s)}-1\Bigr).
   $$}
   \bigskip
   {\bf Proof} \ Indeed, by the definition of $I_{s,s}$, we get
   $$
   \log{I_{s,s}(a,b)\over S(a,b)}=-{1\over s}+{a^s\log a-b^s\log b\over
   a^s-b^s}-{a\log a+b\log b\over a+b}
   $$
   $$
   =-{1\over s}+ab{a^{s-1}+b^{s-1}\over a+b}{\log a-\log b\over
   a^s-b^s}={1\over s}\Bigl({(a^s)^{1-1/s}+(b^s)^{1-1/s}\over (a^s)^{-1/s}+(b^s)^{-1/s}}{\log a^s-\log b^s\over a^s-b^s}-1
   \Bigr)
   $$
   $$
   = {1\over s}\Bigl({l_{-{1\over
s}}(a^s,b^s)\over L(a^s,b^s)}-1\Bigr).
   $$
   \bigskip
   Now, putting $s=3$ in the above identity and applying Lemma
   3.3, the proof follows immediately.
   \bigskip
   Therefore, by Property d. of Proposition 3.1, for $r,s\in [-3,3]$ we get
   $$
   I_{r,s}\le I_{3,3}\le S.
   $$
   \bigskip
   Also, since for fixed $s, s>3$,
   $$
   \sigma(I_{s,s})=2e^{-1/s}<2=\sigma(S),
   $$
   it follows by Theorem 2.0 that the mean $S$ is the right
   cancelling mean for $\{I_{s,s}\}$.
   \bigskip
   Similarly,
   $$
   I_{r,s}\ge I_{-3,-3},
   $$

   and the left hand side of Theorem 3.2 follows from
   easy-checkable relations
   $$
   I_{-s,-s}(a,b)={ab\over I_{s,s}(a,b)}, \ \ a^{b\over a+b}b^{a\over
   a+b}={ab\over S(a,b)}.
   $$

\vskip 1cm

 {\bf 4. Discussion and some open questions}
 \bigskip
Obviously, the right cancelling mean $S^*(\Delta)$ (respectively,
the left cancelling mean $S_*(\Delta)$) is not unique. For
instance, $T(a,b)={1\over2}(S^*(\Delta)+\max (a,b)), \ T\in\Omega$
is also cancelling mean for the class $\Delta$.
\bigskip
Therefore, the mean $S$ is not an exclusive right cancelling mean
in the above assertions. Moreover, we can construct a whole class
of means which may replace the mean $S$ as the right cancelling
mean.
\bigskip

 {\bf Theorem 4.1} \ {\it For $r>-1$, each term of the family of means $K_r$,
 $$
 K_r=K_r(a,b):=\Bigl({a^{r+1}+b^{r+1}\over a+b} \Bigr)^{1/r},
 \ K_0=S,
 $$

 can be taken as the right cancelling mean for the class
 $\{M_s\}$.
 \bigskip
 {\bf Proof} \ We shall prove first that $K_r$ is monotone
 increasing in $r\in\Bbb R$. For this aim, consider the weighted
 arithmetic mean $A_{p,q}(x,y):=px+qy$, where $p,q$ are arbitrary positive numbers such that
 $p+q=1$. Since
 $$
[\log A_{p,q}(x,y)]_{xy}=-{pq\over (px+qy)^2},
 $$

 we conclude that
 $$
\tilde{A}_r(p,q;a,b):=(pa^r+qb^r)^{1/r},
 $$

 is monotone increasing in $r\in\Bbb R$.
 \bigskip
 Hence, the relation
 $$
 \tilde{A}_r({a\over a+b},{b\over a+b};a,b)=K_r(a,b),
 $$

 yields the proof.
 \bigskip
 Now, since for fixed $r>-1$,
 $$
M_0=G\le A=K_{-1}\le K_r,
 $$

and $\sigma (K_r)=2$, it follows that $K_r$ is the right
cancelling mean for the class $\{M_s\}$ analogously to the proof
of Theorem 2.1.
\bigskip
Finally, we propose two open questions concerning the above
matter.
\bigskip
{\bf Q1} \ {\it Does there exists $\min (S^*(A_s))$?}
\bigskip
Denote by $\{K'_r\}$ the subset of $\{K_r\}$ with $r>-1$ i.e.
$\sigma (K'_r)=2$. Then $\max (S_*(K'_r))=K_{-1}=A$.
\bigskip
{\bf Q2} \ {\it Does there exists a non-trivial right cancelling
mean for the class $\{K'_r\}$ ?}

 \vskip 1cm

{\bf References}
\bigskip
[1] \ Hardy, G.H., Littlewood, J. E., P\"{o}lya, G.: {\it
Inequalities}, Camb. Univ. Press, Cambridge (1978).
\bigskip
[2] \ B. C. Carlson, {\it The logarithmic mean}, Amer. Math.
Monthly, 79 (1972), pp. 615-618.
\bigskip
[3] \ P.A. Hasto, {\it Optimal inequalities between Seiffert's
mean and power means}, Math. Ineq. Appl. Vol. 7, No. 1 (2004), pp.
47-53.
\bigskip

 [4] \ T. P. Lin, {\it The power mean and the logarithmic
mean}, Amer. Math. Monthly, 81 (1974), pp. 879-883.
\bigskip
[5] \  Stolarsky, K.: \ {\it Generalizations of the logarithmic
mean}, Math. Mag. 48 (1975), pp. 87-92.
\bigskip
[6] \ Simic, S. : {\it On logarithmic convexity for differences of
power means}, J. Ineq. Appl. Article ID 37359 (2007).
\bigskip
[7] \ Simic, S. : {\it An extension of Stolarsky means to the
multivariable case}, Inter. J. Math. Math. Sci. Article ID 432857
(2009).
\bigskip
[8] \ Simic, S. : {\it On some intermediate mean values}, Inter.
J. Math. Math. Sci. Article ID 283127 (2012),
doi:10.1155/2012/283127.
\bigskip
[9] \ Yang, Z-H. : {\it On the homogeneous functions with two
parameters and its monotonicity}, J. Inequal. Pure and Appl.
Math., Vol. 6, Issue 4 (2005).
\bigskip
[10] \ Yang, Z-H. : {\it Sharp bounds for the second Seiffert mean
in terms of power mean}, arXiv: 1206.5494v1 [math. CA] (2012).
\bigskip
[11] \ Qiu, Y-F., Wang, M-K., Chu, Y-M. and Wang, G. : {\it Two
sharp bounds for Lehmer mean, identric mean and logarithmic mean
}, J. Math. Inequal., Vol. 5, No. 3 (2011), pp. 301-306.

\end